\documentclass[a4paper]{article}
\usepackage[utf8]{inputenc}  
\usepackage[T1]{fontenc}

\usepackage{lmodern} \normalfont 
\DeclareFontShape{T1}{lmr}{bx}{sc} { <-> ssub * cmr/bx/sc }{}

\usepackage{amsthm}
\usepackage{amsmath}
\usepackage{amssymb}
\usepackage{authblk}

\usepackage{enumitem}

\setitemize[0]{label=$\bullet$}
\setitemize[1]{label=$\circ$}

\usepackage[nottoc, notlof, notlot]{tocbibind}

\usepackage{pgf,tikz}
\usetikzlibrary{arrows}
\usetikzlibrary{patterns}

\usepackage[all]{xy}

\usepackage{subcaption}
\usepackage{geometry}
\usepackage{ifthen}

\usepackage[colorlinks=true,linkcolor=blue,citecolor=green]{hyperref}

\newcommand{\R}{\mathbb{R}}

\newcommand{\Q}{\mathbb{Q}}
\newcommand{\C}{\mathbb{C}}

\renewcommand{\P}{\mathbb{P}}
\renewcommand{\O}{{\cal O}}

\def\bsys{\left\{\begin{array}}
\def\esys{\end{array}\right.}

\newcommand{\eps}{\varepsilon}

\theoremstyle{plain}
\newtheorem{theorem}{Theorem}[section]
\newtheorem{lemma}[theorem]{Lemma}
\newtheorem{proposition}[theorem]{Proposition}
\newtheorem{corollary}[theorem]{Corollary}

\theoremstyle{definition}
\newtheorem{definition}[theorem]{Definition}

\newtheorem{remark}[theorem]{Remark}
\newtheorem{example}[theorem]{Example}
\newtheorem{nexample}[theorem]{Non-exemple}

\newcommand{\lem}[2][None]
{\begin{lemma}\ifthenelse{\equal{#1}{None}}{}{\label{#1}}
#2
\end{lemma}}

\newcommand{\prop}[2][None]
{\begin{proposition}\ifthenelse{\equal{#1}{None}}{}{\label{#1}}
#2
\end{proposition}}

\newcommand{\cor}[2][None]
{\begin{corollary}\ifthenelse{\equal{#1}{None}}{}{\label{#1}}
#2
\end{corollary}}

\newcommand{\theo}[2][None]
{\begin{theorem}\ifthenelse{\equal{#1}{None}}{}{\label{#1}}
#2
\end{theorem}}

\newcommand{\defi}[2][None]
{\begin{definition}\ifthenelse{\equal{#1}{None}}{}{\label{#1}}
#2
\end{definition}}

\newcommand{\rem}[2][None]
{\begin{remark}\ifthenelse{\equal{#1}{None}}{}{\label{#1}}
#2
\end{remark}}

\newcommand{\expl}[2][None]
{\begin{example}\ifthenelse{\equal{#1}{None}}{}{\label{#1}}
#2
\end{example}}

\newcommand{\nexpl}[2][None]
{\begin{nexample}\ifthenelse{\equal{#1}{None}}{}{\label{#1}}
#2
\end{nexample}}

\newcommand{\demo}[1]{\begin{proof} #1
\end{proof}}

\newcommand{\demode}[2]{\begin{proof}[Proof of #1] #2
\end{proof}}

\newcounter{exerc}

\newcounter{quest}

\newcommand{\qst}[2][None]{\refstepcounter{quest}\ifthenelse{\equal{#1}{None}}{}{\label{#1}} \smallskip \noindent {\bf Question \thequest\;} {#2} \smallskip}

\newcommand{\exo}[2][None]{\refstepcounter{exerc}\ifthenelse{\equal{#1}{None}}{}{\label{#1}} \smallskip \noindent {\bf Exercise \theexerc\;} {#2} \smallskip}

\newcommand{\eqnum}{\underset{\mbox{{\footnotesize num}}}{\equiv}}

\author{Cécile Gachet\footnote{Université Côte d'Azur, CNRS, LJAD, France; gachet@unice.fr}}
\title{Positivity of higher exterior powers of the tangent bundle} 

\begin{document}

\maketitle

\begin{abstract}
We prove that a smooth projective variety $X$ of dimension $n$ with strictly nef third, fourth or $(n-1)$-th exterior power of the tangent bundle is a Fano variety. Moreover, in the first two cases, we provide a classification for $X$ under the assumption that $\rho(X)\ne 1$.
\end{abstract}

\section{Introduction}

Positivity notions are numerous in algebraic geometry: a line bundle can be considered positive, e.g., if it is very ample, ample, strictly nef, nef, big, semiample, effective, pseudoeffective... Some of these notions relate: a very ample line bundle is ample, an ample line bundle is strictly nef and big, a strictly nef line bundle (i.e., a line bundle that has positive intersection with any curve) is nef, a nef line bundle and an effective line bundle are pseudoeffective. These positivity notions, as they tremendously matter in algebraic geometry, have been the subject of a lot of work, to which the books by Lazarsfeld \cite{LazBis,Laz} are a great introduction.
Proving new relationships between these various positivity notions is however a rather naive ambition, if not under strong additional assumptions.

From this perspective, the conjecture by Campana and Peternell \cite{CP91} is surprising: they predict that, if $X$ is a smooth projective variety, and the anticanonical divisor $-K_X$ is strictly nef, then $-K_X$ is ample, i.e., $X$ is a Fano manifold. Their conjecture was in fact proven in dimension 2 and 3, by Maeda and Serrano \cite{Maeda, Serra}. As all Fano manifolds are rationally connected \cite{Camp,KMM}, an interesting update on the conjecture is the recent proof by Li, Ou and Yang \cite[Theorem 1.2]{Main} that if $X$ is a smooth projective variety, and the anticanonical divisor $-K_X$ is strictly nef, then $X$ is rationally connected. Their proof uses important results on the Albanese map of varieties with nef anticanonical bundle. Such varieties have been extensively studied too \cite{DPS,Zhang,PeternellSerrano,Demailly,CaoHoe1,Cao,CaoHoe2}.

Positivity notions extend to vector bundles \cite[Definition 6.1.1]{Laz} in the following fashion: a vector bundle $E$ is stricly nef if the associated line bundle $\O_{\P(E)}(1)$ is strictly nef on $\P(E)$. Instead of asking about the positivity of the top exterior power of the tangent bundle, $-K_X=\bigwedge^{\dim(X)}T_X$, it makes sense to ask about the positivity of intermediate exterior powers $\bigwedge^r T_X$, for $1\le r\le \dim(X)-1$. 

For $r=1$, it is known since Mori \cite{Mor1} that projective spaces are the only smooth projective varieties with ample tangent bundle. They are also the only smooth projective varieties with strictly nef tangent bundle, by \cite[Theorem 1.4]{Main}. Varieties with nef tangent bundle are, on the other hand, governed by another conjecture of Campana and Peternell \cite{CP91} which has received a lot of attention: see the survey \cite{SurveyNef}, and {\it inter alia} \cite{CP91,DPS,WataNef,Kanemitsu1,Kanemitsu2,DynkinNef,YangToric,LiHorospherical,Demailly17,WataCoindex,WataPosChar}.

For $r=2$, it has been proven that varieties with ample second exterior power of the tangent bundle are projective spaces and quadric hypersurfaces \cite{ChoSato}, varieties with strictly nef second exterior power of the tangent bundle alike.

\theo[1.1]{\textup{\cite[Theorem 1.5]{Main}} Let $X$ be a smooth projective variety of dimension $n\ge 2$, such that $\bigwedge^2 T_X$ is strictly nef. Then $X$ is isomorphic to the projective space $\P^n$, or to a smooth quadric hypersurface $Q^n$.}

Partial results were obtained under the nef assumption \cite{Wata,Schmitz}. 

These results lead us to the following questions.

\medskip

\noindent{\bf Question 1.} Let $X$ be a smooth projective variety of dimension $n$. Suppose that $\bigwedge^r T_X$ is strictly nef for some integer $1\le r\le n$. Is $X$ a Fano variety?

\medskip

\noindent{\bf Question 2.} Let $X$ be a smooth projective variety of dimension $n$. Suppose that $\bigwedge^r T_X$ is nef for some integer $1\le r< n$, and that $X$ is rationally connected. Is $X$ a Fano variety?

\medskip

Note that an affirmative answer to the second question would imply an affirmative answer to the first question, by \cite[Theorem 1.2]{Main}. Also note that the second question is answered negatively for $r=n$, as there are smooth rationally connected threefolds with $-K_X$ nef but not semiample \cite{Zhixin}. The first question is answered affirmatively for smooth toric varieties by \cite{Schmitz}.
In this paper, we answer the second question for $r=n-1$.

\theo[theo-n-1main]{Let $X$ be a smooth projective variety of dimension $n\ge 2$ such that the vector bundle $\bigwedge^{n-1} T_X$ is nef and $X$ is rationally connected. Then $X$ is a Fano variety.}

This theorem is reminiscent of \cite[Proposition 3.10]{DPS}, which states a dichotomy for varieties $X$ with nef tangent bundle: either $X$ is a Fano manifold, or $\chi(X,\O_X)=0$. The proof similarly involves Chern classes inequalities and the Hirzebruch-Riemann-Roch formula.
Note that, building on this theorem, \cite[Proposition 1.4]{WataNew} very recently gave an affirmative answer to Question 2 in general. 

Theorem \ref{1.1} is based on the results of \cite{CMSB} and \cite{DeHoe}, which instead of the assumption on $\bigwedge^2 T_X$, feature a much weaker assumption on the length of rational curves. In a similar spirit, we provide the following partial characterizations and their corollaries. 

\theo[theo-3mainnew]{Let $X$ be a smooth projective rationally connected variety of dimension $n\ge 4$ such that for each rational curve $C$ in $X$, we have $-K_X\cdot C\ge n-1$. Then either $X\simeq \P^2\times\P^2$, or $X$ is a Fano variety of Picard rank $\rho(X) = 1$.}

\cor[cor-3main]{Let $X$ be a smooth projective variety of dimension at least $4$ such that the vector bundle $\bigwedge^3 T_X$ is strictly nef. Then either $X\simeq \P^2\times\P^2$, or $X$ is a Fano variety of Picard rank $\rho(X) = 1$.}

Let us briefly discuss the case when $\rho(X)=1$. We know that, if $X$ is a cubic or a complete intersection of two quadrics in $\P^n$, the vector bundle $\bigwedge^3 T_X$ is ample. These are two examples of del Pezzo manifolds, i.e. Fano $n$-folds of Picard rank 1 and of index $n-1$. However, we do not know whether other del Pezzo manifolds have strictly nef $\bigwedge^3 T_X$, or whether varieties with strictly nef $\bigwedge^3 T_X$ are in general del Pezzo manifolds. We can hardly hope for a characterization of Fano manifolds of Picard rank one on which $-K_X\cdot C \ge n-1$ for every rational curve $C$, and it is moreover not clear how to use the positivity of $\bigwedge^3 T_X$ beyond that length inequality, {\it cf.} Lemma \ref{degrantican}.

\theo[theo-4mainnew]{Let $X$ be a smooth projective rationally connected variety of dimension $n\ge 6$ such that for each rational curve $C$ in $X$, we have $-K_X\cdot C\ge n-2$. Then either $X$ is isomorphic to $\P^3\times\P^3$
or $X$ is a Fano variety of Picard rank $\rho(X) = 1$.}

Studying the possibilities in dimension 5 by hand yields the following result.

\cor[cor-4main]{Let $X$ be a smooth projective variety of dimension at least $5$ such that the vector bundle $\bigwedge^4 T_X$ is strictly nef. Then either $X$ is isomorphic to one of the following Fano varieties
\begin{center}
 $\P^2\times Q^3$;\;
 $\P^2\times\P^3$;\;
 $\P(T_{\P^3})$;\;
 $\mathrm{Bl}_{\ell}(\P^5)=\P(\O_{\P^3}\oplus \O_{\P^3}\oplus \O_{\P^3}(1))$;
 $\P^3\times\P^3$
\end{center}
or $X$ is a Fano variety of Picard rank $\rho(X) = 1$.}

These two corollaries were to our knowledge unknown even under the stronger, more classical assumption that $\bigwedge^3 T_X$ or $\bigwedge^4 T_X$ is ample. The proof of both theorems goes by classifying possible Mori contractions for $X$. A delicate point is that, while we know that our varieties $X$ with $\rho(X)\ge 2$ admit one Mori contraction by the Cone Theorem, we need to construct by hand a second Mori contraction, e.g., to control higher-dimensional fibres in case of a first fibred Mori contraction. Depending on circumstances, we use unsplit covering families of deformations of rational curves, and a result by Bonavero, Casagrande and Druel \cite{BonCaD}, or, if $X$ has the right dimension, Theorem \ref{theo-n-1main}, to produce this second Mori contraction.

\bigskip

\noindent{\bf Acknowledgments.} I am grateful to my advisor A. Höring for regular helpful discussions, to S. Tanimoto for pointing out that the complete intersection of two quadrics in a projective space should satisfy Corollary \ref{cor-3main}, and to J. Cao for suggesting the second question in the introduction.

\bigskip

\noindent{\bf Conventions.}
We work over the field of complex numbers $\C$. Varieties (and in particular curves) are always assumed irreducible and reduced. We use the expressions ``smooth projective variety'' and ``projective manifold'' interchangeably. We refer to \cite{Debarre} for birational geometry, in particular Mori theory, \cite{LazBis,Laz} for positivity notions, \cite{Koll} for rational curves and their deformations. We write $c_i(X)=c_i(T_X)$ for the Chern classes of the tangent bundle of $X$.

\section{A first lemma}

We start with a simple lemma.

\lem[degrantican]{Let $X$ be a smooth projective variety of dimension $n$ such that $\bigwedge^r T_X$ is strictly nef, for some $1\le r\le n-1$. Then any rational curve $C$ in $X$ satisfies 
$$-K_X\cdot C\ge n + 2 - r.$$}

\demo{The proof goes as \cite[Proof of Theorem 1.5]{Main}. Let $f:\P^1\to C$ be the normalization of the curve. Write
$$f^*T_X\simeq\O_{\P^1}(a_1)\oplus\ldots\oplus\O_{\P^1}(a_n),$$

\noindent with $(a_i)_{1\le i\le n}$ ordered increasingly. It holds $a_n\ge 2$, as $T_{\P^1}$ maps non-trivially to $f^*T_X$, and we have $a_1+\ldots+a_r> 0$ because $\O_{\P^1}(a_1+\ldots+a_r)$ is a direct summand of the strictly nef vector bundle $\bigwedge^r f^*T_X$. In particular, $a_{r+1} \ge a_r\ge 1$. Hence,

$$-K_X\cdot C=\deg f^*(-K_X)=a_1+\ldots+a_n\ge 1 + n-r-1 + 2 = n+2-r.$$
}

This result is all the more valuable as, by \cite[Theorem 1.2]{Main}, if $X$ is a smooth projective variety of dimension $n$ such that $\bigwedge^r T_X$ is strictly nef, then it is rationally connected, in particular, it contains numerous rational curves.

We will also need the following result.

\lem[degranticanbis]{Let $X$ be a smooth projective variety of dimension $n$ such that $\bigwedge^r T_X$ is nef, for some $1\le r\le n-1$. Then any rational curve $C$ in $X$ satisfies 
$-K_X\cdot C\ge 2.$}

\demo{Let $f:\P^1\to C$ be the normalization of the curve. Write
$$f^*T_X\simeq\O_{\P^1}(a_1)\oplus\ldots\oplus\O_{\P^1}(a_n),$$

\noindent with $(a_i)_{1\le i\le n}$ ordered increasingly. It holds $a_n\ge 2$, as $T_{\P^1}$ maps non-trivially to $f^*T_X$, and we have $a_1+\ldots+a_r\ge 0$ because $\O_{\P^1}(a_1+\ldots+a_r)$ is a direct summand of the nef vector bundle $\bigwedge^r f^*T_X$. Hence, $a_{r+1}\ge a_r\ge 0$, and summing up those inequalities, we obtain the estimate 
$$-K_X\cdot C = a_1+\ldots+a_n\ge 2.$$}

 \section{Results on $\bigwedge^{n-1} T_X$}\label{sec-n-1}

The following lemma is the main step in the proof of Theorem \ref{theo-n-1main}.

 \lem[lem-1]{Let $X$ be a projective $n$-dimensional manifold such that $\bigwedge^{n-1} T_X$ is nef and $X$ is rationally connected. Then $-K_X$ is nef and big.}
 
 \demo{By \cite[Theorem 6.2.12(iv)]{Laz}, the anticanonical divisor $-K_X$ is nef. By the Hirzebruch-Riemann-Roch formula, there is a homogeneous polynomial $P$ of degree $n$ in $\Q[X_1,\ldots,X_n]$ with grading $\deg\, X_i = i$ such that 
 $$\chi(X,\O_X)=P(c_1(X),\ldots,c_n(X)).$$

Note that, as $\bigwedge^{n-1} T_X = \Omega_X^1\otimes\O_X(-K_X)$, and by \cite[Remark 3.2.3(b)]{Fulton}, we have

\begin{equation}\tag{$\ast$}\label{eq-chern}
c_i\left(\bigwedge^{n-1} T_X\right) = \sum_{j=0}^i (-1)^j {n-j \choose i-j}
c_j(X)c_1(-K_X)^{i-j}.
\end{equation}

Let us show by induction that $c_i(X)$ is a rational polynomial in the $c_j(\bigwedge^{n-1} T_X)$, for $0\le j\le i$. Indeed, $c_1(X) = \frac{1}{n}c_1(\bigwedge^{n-1} T_X)$. Assume now that for some $i$, for all $0\le j\le i$, there is a polynomial $P_j\in\Q[X_1,\ldots,X_j]$ such that $c_j(X)=P_j(c_1(\bigwedge^{n-1} T_X),\ldots,c_j(\bigwedge^{n-1} T_X))$. 
Then, setting 
$$P_{i+1}(X_1,\ldots,X_{i+1})=(-1)^{i+1} X_{i+1}
- \sum_{j=0}^i (-1)^{i+j+1} {n-j \choose i+1-j} P_j(X_1,\ldots,X_j)(P_1(X_1))^{i+1-j},$$ 
we have
$c_{i+1}(X)= P_{i+1}(c_1(\bigwedge^{n-1} T_X),\ldots,c_{i+1}(\bigwedge^{n-1} T_X))$ by (\ref{eq-chern}). This perpetuates the induction.

In particular, we have
 $$\chi(X,\O_X)=P\left(P_1\left(c_1\left(\bigwedge^{n-1} T_X\right)\right),\ldots,P_n\left(c_1\left(\bigwedge^{n-1} T_X\right),\ldots,c_n\left(\bigwedge^{n-1} T_X\right)\right)\right),$$
which is a homogeneous polynomial of degree $n$ in $c_1(\bigwedge^{n-1} T_X),\ldots,c_n(\bigwedge^{n-1} T_X)$.

Now, if we suppose that $-K_X$ is not big, then $c_1(\bigwedge^{n-1}T_X)$ is not big. Thus, \cite[Corollary 2.7]{DPS} implies $\chi(X,\O_X)=0$. But on the other hand, $X$ is rationally connected, so $\chi(X,\O_X)=1$, a contradiction.}

\rem{If $n=4$, we cannot write $c_3(X)$ as a polynomial in 
\begin{align*}
c_1\left(\bigwedge^{n-2} T_X\right) &=3c_1(X), \\
c_2\left(\bigwedge^{n-2} T_X\right) &=3c_1(X)^2+2c_2(X), \\
c_3\left(\bigwedge^{n-2} T_X\right) &=c_1(X)^3+4c_1(X)c_2(X),
\end{align*} these formulas coming from \cite[4.5.2]{Iena}.}

\lem[lem-ample]{Let $X$ be a projective $n$-dimensional manifold such that $\bigwedge^{n-1} T_X$ is nef and $X$ is rationally connected. Then $-K_X$ is ample.}

\demode{Theorem \ref{theo-n-1main}}{By Lemma \ref{lem-1}, $-K_X$ is nef and big. By the base-point-free theorem \cite[Theorem 7.32]{Debarre}, we dispose of an integer $m$ such that $-mK_X$ is globally generated. Let $\eps:X\to Z$ be the $|-mK_X|$-morphism. 

Suppose that it is not finite. By \cite[Theorem 2]{KawaLength}, any irreducible component $E$ of the exceptional locus is covered by rational curves that are contracted by $\eps$. Let $C$ be one of them: we have $0=-K_X\cdot C\ge 2$ by Lemma \ref{degranticanbis}, a contradiction. So $-K_X$ is ample.}

\section{Studying Mori contractions}\label{sec-general}

The strategy for proving Theorems \ref{theo-3mainnew} and \ref{theo-4mainnew} is to show that there are only few possible Mori contractions for $X$. 
In the following, if $R$ is an extremal ray of the Mori cone $\overline{NE}(X)$, its {\it length} denoted by $\ell(R)$ is defined to be the minimal value of $-K_X\cdot C$, for a rational curve $C$ with class in $R$. A Mori contraction is said to be {\it of length} $\ell$ if it is a contraction of an extremal ray $R$ with $\ell(R)=\ell$. Note that here and throughout this paper, we exclusively work with elementary Mori contractions, and simply refer to them as Mori contractions.

\subsection{Small contractions}

\lem[paspet]{Let $1\le r\le 4$. Let $X$ be a smooth projective variety of dimension at least $r + 1$ such that $\bigwedge^r T_X$ is strictly nef. Then $X$ has no small contraction.}

\demo{Let $n$ be the dimension of $X$. Let $\varphi : X\to Y$ be a birational contraction, $E$ be an irreducible component of the exceptional locus, $F$ an irreducible component of the general fibre of $\varphi|_E$, and $R$ the corresponding extremal ray. 
Applying Ionescu-Wi\'snewski inequality \cite[Theorem 0.4]{Io}, \cite[Theorem 1.1]{Wisn} together with Lemma \ref{degrantican} yields
$$\dim E + \dim F\ge n +\ell(R) - 1 \ge 2n+1-r.$$

\noindent Since $r\le 4$, we have $\dim E\ge n-1$, and thus $\varphi$ is a divisorial contraction.}

\subsection{Fibred Mori contractions}

We move on to studying fibred Mori contractions.

\subsubsection{Generalities about fibred Mori contractions}

If $X$ is a normal projective variety, and $C$ is a rational curve in $X$, we may denote by ${\cal V}$ its family of deformations, that is an irreducible component of $\mathrm{Chow}(X)$ containing the point corresponding to $C$. Denoting by $\phi:\mathrm{Univ}(X)\to\mathrm{Chow}(X)$ the universal family and by $\mathrm{ev}:\mathrm{Univ}(X)\to X$ the evaluation map, we define
$$\mathrm{Locus}({\cal V}):=\mathrm{ev}(\phi^{-1}({\cal V}))\subset X.$$
We say that ${\cal V}$ is {\it covering} if $\mathrm{Locus}({\cal V})=X$.

\noindent
We say that ${\cal V}$ is {\it unsplit} if it only parametrizes irreducible cycles.

For $x\in\mathrm{Locus}({\cal V})$, we define ${\cal V}_x:=\phi(\mathrm{ev}^{-1}(x))$ the family of deformations of $C$ through $x$. We finally define
$$\mathrm{Locus}({\cal V}_x):=\mathrm{ev}(\phi^{-1}({\cal V}_x))\subset X.$$

We use families of deformations of rational curves to prove the following proposition.

\prop[fibpet]{Let $X$ be a smooth projective rationally connected variety of dimension $n$. Let $r$ be an integer with $1\le r\le n-1$. Suppose that $-K_X\cdot C\ge n+2-r$ for any rational curve $C$ in $X$. Suppose that there is a fibred Mori contraction $\pi:X\to Y$ with $\dim Y > 0$. Then the general fibre has dimension at most $r-1$. 

If equality holds, then there is a rational curve $C$ in $X$, not contracted by $\pi$, whose family of deformations ${\cal V}$ is unsplit covering and satisfies $\dim\mathrm{Locus}({\cal V}_x)=n+1-r$ for $x\in\mathrm{Locus}({\cal V})$ general.}

The proof relies on the following lemmas.

\lem[lem-finite]{Let $X$ be a smooth projective variety. Suppose that $X$ has a fibred Mori contraction $\pi:X\to Y$ with $\dim Y > 0$, and let $C$ be a rational curve such that $\pi(C)\ne \{\mathrm{pt}\}$ and such that its family of deformations ${\cal V}$ is unsplit. Then, for any $x\in\mathrm{Locus}({\cal V})$, 
$$\dim\mathrm{Locus}({\cal V}_x)\le\dim Y.$$}

\demode{Lemma \ref{lem-finite}}{We claim that $\pi|_{\mathrm{Locus}({\cal V}_x)}$ is finite onto its image. If it is not, it contracts a curve $B$ to a point: for some ample divisor $H$ on $Y$, we have $B\cdot \pi^* H=0$.
By \cite[Lemma 4.1]{ACOMuk}, the numerical class of $B\subset\mathrm{Locus}({\cal V}_x)$ is a multiple of $[C]\in N_1(X)_{\Q}$, whence $C\cdot\pi^*H=0$, which is a contradiction.
So $\pi|_{\mathrm{Locus}({\cal V}_x)}$ is finite onto its image: this implies $\dim \mathrm{Locus}({\cal V}_x)\le \dim Y $.}

\lem[lem-unsplit]{Let $X$ be a smooth projective variety. Suppose that $-K_X\cdot C > 0$ for every rational curve $C\subset X$. Suppose that $X$ has a fibred Mori contraction $\pi:X\to Y$ with $\dim Y > 0$, and let $C$ be a rational curve such that $\pi(C)\ne \{\mathrm{pt}\}$ and such that 
$$-K_X\cdot C=\min\{-K_X\cdot B\mid B\mbox{ rational curve in }X, \pi(B)\ne\{\mathrm{pt}\}\}.$$
Then the family of deformations of $C$ is unsplit.}

\demode{Lemma \ref{lem-unsplit}}{Let ${\cal V}$ be the family of deformations of $C$. 
Suppose that it is splitting. By \cite[Explanation IV.2.2]{Koll}, we have
$$C\eqnum\sum_i a_iC_i,$$
with rational curves $C_i$ and coefficients $a_i\ge 1$ such that $\sum_{i} a_i\ge 2$.
Since $-K_X$ is positive on rational curves, we have $-K_X\cdot C_i < -K_X\cdot C$ for all $i$. So, by minimality of $-K_X\cdot C$, the fibration $\pi$ contracts all curves $C_i$. Let $H$ be an ample divisor on $Y$. We obtain $\sum_i a_iC_i\cdot \pi^*H=0$, a contradiction.}

\demode{Proposition \ref{fibpet}}{Since $X$ is rationally connected and $-K_X$ is Cartier, we dispose of a rational curve $C$ such that $\pi(C)\ne\{\mathrm{pt}\}$ and $-K_X\cdot C\ge n+2-r\ge 3$ is minimal with this condition.
Let ${\cal V}$ be the corresponding family of deformations. By Lemma \ref{lem-unsplit}, it is unsplit.

Fix $x\in \mathrm{Locus}({\cal V})$ general. 
 By \cite[Proposition IV.2.6]{Koll} and our assumption, we derive
$$\dim \mathrm{Locus}({\cal V})+ \dim\mathrm{Locus}({\cal V}_x)
\ge -K_X\cdot C + n -1
\ge 2n +1 -r.$$
So $\dim \mathrm{Locus}({\cal V}_x) \ge n + 1 -r.$

Let $d$ denote the dimension of the general fibre of $\pi$. Then, by Lemma \ref{lem-finite},
$$d\le n-\dim\mathrm{Locus}({\cal V}_x)\le r-1.$$

As for the equality case, if $d=r-1$, then $\dim\mathrm{Locus}({\cal V}_x)= n-r+1$, and so $C$ is such a rational curve as we claimed existed in the equality case of the proposition.}

Proposition \ref{fibpet} has an important consequence.

\cor[cor-evid]{Let $X$ be a smooth projective rationally connected variety of dimension $n$ such that, for some integer $r$ with $1\le r\le n-1$, one has $-K_X\cdot C\ge n+2-r$ for any rational curve $C\subset X$. Suppose that there is a fibred Mori contraction $\pi:X\to Y$ with $\dim Y > 0$. Then $n\le 2r-2$. 

If equality holds, then a general fibre of $\pi$ has dimension $r-1$, and there is a rational curve $C$ in $X$, not contracted by $\pi$, whose family of deformations ${\cal V}$ is unsplit covering and satisfies $\dim\mathrm{Locus}({\cal V}_x)=n+1-r$ for $x\in\mathrm{Locus}({\cal V})$ general.}

\demo{Let $F$ be a general fiber of $\pi$. By Proposition \ref{fibpet}, we have $r-1\ge \dim F$. Adding $n$ to both sides and applying Ionescu-Wi\'snewski inequality (with the exceptional locus $E=X$ of dimension $n$), it holds
$$n+r-1 \ge n + \dim F\ge n +\ell(R) - 1 \ge 2n+1-r.$$

\noindent If there is an equality, then $\dim F = r-1$, and so we are in the equality case of Proposition \ref{fibpet}. In particular, we can find a rational curve $C$ in $X$ that is not contracted by $\pi$, whose family of deformations ${\cal V}$ is unsplit and satisfies $\dim\mathrm{Locus}({\cal V}_x)=n+1-r$ for $x\in\mathrm{Locus}({\cal V})$ general. By \cite[Proposition IV.2.6]{Koll} and Lemma \ref{degrantican} again, we have
$$\dim \mathrm{Locus}({\cal V})
\ge -K_X\cdot C + n -1 - \dim\mathrm{Locus}({\cal V}_x)
\ge 2n +1 -r - n - 1 + r =n,$$
so ${\cal V}$ is indeed a covering family.}

\subsubsection{Fibred Mori contractions for certain varieties of even dimension}

The set-up for this paragraph is the following. Let $r\ge 3$ be an integer. Let $X$ be a smooth projective rationally connected variety of dimension $2r-2$ such that $-K_X\cdot C\ge r$ for any rational curve $C\subset X$. Suppose that there is a fibred Mori contraction $\pi:X\to Y$ with $\dim Y > 0$.
Let us classify what happens.

\lem[prop-p2p2]{Let $r\ge 3$ be an integer. Let $X$ be a smooth projective rationally connected variety of dimension $2r-2$ such that $-K_X\cdot C\ge r$ for any rational curve $C\subset X$. Suppose that there is a fibred Mori contraction $\pi:X\to Y$ with $\dim Y > 0$. Then there is another equidimensional fibred Mori contraction $\varphi:X\to Z$ with $\dim Z = r -1$.}

\demo{We are in the case of equality of Corollary \ref{cor-evid}. In particular, the general fibre $F$ of $\pi$ has dimension $r-1$, and there is a rational curve $C$ in $X$ that is not contracted by $\pi$ whose family of deformations ${\cal V}$ is unsplit covering and satisfies $\dim\mathrm{Locus}({\cal V}_x)= r-1\ge (2r-2)-3=\dim X - 3$.

By \cite[Theorem 2, Proposition 1(i)]{BonCaD}, there is a fibred Mori contraction $\varphi :X\to Z$ whose fibres exactly are the ${\cal V}$-equivalence classes. By the equality case in Corollary \ref{cor-evid}, the general fibre of $\varphi$ has dimension $r-1$.

Let $G$ be a fibre of $\varphi$. We claim that $\pi|_G$ is finite. Indeed, if it is not, then there is a curve $B\subset G$ that is contracted by $\pi$. The curve $B$ lies in a ${\cal V}$-equivalence class, so by \cite[Remark 1]{BonCaD}, as ${\cal V}$ is unsplit, $B$ is numerically equivalent to a multiple of $C$, so it cannot be contracted by $\pi$, a contradiction! So $\pi|_G$ is finite onto its image, which is contained in $Y$, so $\dim G\le \dim Y = r -1$.

So $\varphi$ is indeed equidimensional.}

\prop[prop-genfib]{Let $r\ge 3$ be an integer.  Let $X$ be a smooth projective rationally connected variety of dimension $2r-2$ such that $-K_X\cdot C\ge r$ for any rational curve $C\subset X$. Suppose that there is an equidimensional fibred Mori contraction $\pi:X\to Y$ with $\dim Y = r -1$. Then $X\simeq \P^{r-1}\times\P^{r-1}$.}

This proposition relies on Lemma \ref{fibproj} below.

\defi{Let $X$ and $Y$ be normal projective varieties. We say that a map $\pi : X\to Y$ is a {\it fibration} if it is surjective, has connected fibers, and if we have $0<\dim Y<\dim X$.

Let $\pi:X\to Y$ be a fibration whose general fibre is a projective space. Let $f:\P^1\to C\subset Y$ be a rational curve whose image lies in the smooth locus of $\pi$. The fibre product $\pi_C$ of $\pi$ by $f$ is the projectivization of a bundle $\O_{\P^1}(a_1)\oplus\ldots\oplus\O_{\P^1}(a_k)$, with the $(a_i)$ ordered increasingly. A {\it minimal section} over $C$ is the section $s:\P^1\to X$ of $\pi_C$ corresponding to a quotient $\O_{\P^1}(a_1)$.}

\rem{There may be several minimal sections as soon as $a_1=a_2$.}

\lem[fibproj]{Let $X$ be a smooth projective variety with a fibration $\pi:X\to Y$ whose general fiber is a projective space. Then for any rational curve $f:\P^1\to C\subset Y$ whose image lies in the smooth locus of $\pi$, and for any minimal section $s$ of it, it holds $-K_Y\cdot C\ge -K_X\cdot s(\P^1)$. In particular,

\begin{equation}\tag{$\ast\ast$}\label{eq-min}
-K_Y\cdot C\ge \min\{-K_X\cdot C'\mid C'\mbox{ is a rational curve in }X\}.
\end{equation}

If there is an equality in (\ref{eq-min}), then the base change of $\pi$ by $f$ is isomorphic to $\P({\O_{\P^1}}^{\oplus k})\to \P^1$.

If there is almost an equality, i.e., 
$$-K_Y\cdot C = \min\{-K_X\cdot C'\mid C'\mbox{ is a rational curve in }X\} + 1,$$ 
then the base change of $\pi$ by $f$ is isomorphic to to $\P({\O_{\P^1}}^{\oplus k})\to \P^1$ or to $\P({\O_{\P^1}}^{\oplus k-1}\oplus\O_{\P^1}(1))\to \P^1$.
}

\demo{By Tsen's theorem, the base change $\pi_C$ of $\pi$ by $f$ is the natural projection morphism of the projectivization of a vector bundle $V$ on $\P^1$. 
We write $V\simeq\O_{\P^1}(a_1)\oplus\ldots\oplus\O_{\P^1}(a_k),$ with $(a_i)$ ordered increasingly, and consider $s$ the section of $\pi_C$ satisfying $s^*\O_{\P(V)}(1)=\O_{\P^1}(a_1)$.
The degree of
$\det (s^*\O_{\P(V)}(1)\otimes V^*)$ is non-positive,
equals zero if and only if
$V\simeq \O_{\P^1}(a_1)^{\oplus k}$, and equals $-1$ if and only if $V\simeq {\O_{\P^1}(a_1)}^{\oplus k-1}\oplus\O_{\P^1}(a_1+1)$.

Pulling-back the Euler exact sequence of $\pi_C$ by $s$, we get
$$0\to\O_{\P^1}\to s^*\O_{\P(V)}(1)\otimes V^* \to s^*T_{X/Y}\to 0.$$
Thus, $s^*T_{X/Y}$ has non-positive degree.
We also have the tangent bundle exact sequence:
$$0\to s^*T_{X/Y}\to s^* T_X\to f^*T_Y\to 0,$$
Since $s^*T_{X/Y}$ has non-positive degree, we obtain
$$-K_Y\cdot C\ge -K_X\cdot s(C)\ge \min\{-K_X\cdot C'\mid C'\mbox{ is a rational curve in }X\}.$$

Moreover, if there is an equality, then we have $-K_Y\cdot C = -K_X\cdot s(C)$, and so $V\simeq \O_{\P^1}(a_1)^{\oplus k}$.

If there is almost an equality, then $-K_Y\cdot C = -K_X\cdot s(C)$ or $-K_Y\cdot C = -K_X\cdot s(C)+1$, so $V\simeq \O_{\P^1}(a_1)^{\oplus k}$ or $V\simeq {\O_{\P^1}(a_1)}^{\oplus k-1}\oplus\O_{\P^1}(a_1+1)$.}

\demode{Proposition \ref{prop-genfib}.}{By \cite[Theorem 1.3]{HoeNov}, as $\pi: X \to Y$ is an equidimensional fibration with fibres of dimension $r-1$, and as it is a Mori contraction of length at least $r$ as well, it is a $\P^{r-1}$-bundle.
Let us show that $Y$ is isomorphic to $\P^{r-1}$.
Since $X$ is smooth and a projective bundle over $Y$, the variety $Y$ is smooth.
By Lemma \ref{fibproj}, any rational curve $C$ in $Y$ satisfies $-K_Y\cdot C\ge r.$
Moreover, $X$ is rationally connected, so $Y$ is too. By \cite[Cor.0.4, $1\Leftrightarrow 10$]{CMSB}, we get $Y\simeq\P^{r-1}$.

As $\P^{r-1}$ has trivial Brauer group, there is a vector bundle $V$ of rank $r$ on $Y$ such that $\pi$ identifies with the natural projection $\P(V)\to \P^{r-1}$. Without loss of generality, we can twist $V$ by a line bundle so that $0\le \deg_{\Delta} V|_{\Delta}\le r-1$, for any line $\Delta$ in $\P^{r-1}$.
Let $\Delta$ be a line in $\P^{r-1}$. Then $-K_{\P^{r-1}}\cdot\Delta = r$. By the equality case in Lemma \ref{fibproj}, the restriction $V|_{\Delta}$ is isomorphic to $L^{\oplus r}$ for some line bundle $L$ on $\Delta$. Hence $\deg L=0$, so $L=\O_{\Delta}$. By \cite[Theorem 3.2.1]{Oko}, the vector bundle $V$ is globally trivial. Hence, $X\simeq \P^{r-1}\times\P^{r-1}$.
}

\subsubsection{Fibred Mori contractions for certain fivefolds}

The goal in this section is prove the following result.

\prop[prop-summary]{Let $X$ be a smooth projective fivefold such that $\bigwedge^4 T_X$ is strictly nef. Suppose that $\rho(X)>1$, and that $X$ admits a fibred Mori contraction. Then $X$ is isomorphic to one of the following projective manifolds 
\begin{center}
 $\P^2\times Q^3$;
 $\P^2\times\P^3$;
 $\P(T_{\P^3})$;
 $\P(\O_{\P^3}\oplus \O_{\P^3}\oplus \O_{\P^3}(1))$.
\end{center}
}

We first establish this classification under the simplifying assumption that $X$ has a $\P^2$-bundle structure, instead of a fibred Mori contraction.

\lem[lem-p2bundle]{Let $X$ be a smooth projective rationally connected fivefold and such that, for any rational curve $C\subset X$, one has $-K_X\cdot C \ge 3$. Suppose that $p:X\to Y$ is a $\P^2$-bundle. Then $Y$ is a smooth projective variety, and $X$ is isomorphic to one of the following projective manifolds
\begin{center}
 $\P^2\times Q^3$;
 $\P^2\times\P^3$;
 $\P(T_{\P^3})$;
 $\P(\O_{\P^3}\oplus \O_{\P^3}\oplus \O_{\P^3}(1))$.
\end{center}
}

Among other things, the proof uses the following lemma.

\lem[lem-trivquadric]{Let $V$ be a vector bundle on a smooth quadric hypersurface $Q^n$. If $V$ is trivial on all lines in $Q^n$, then $V$ is trivial.}

\demo{Note that by \cite[Theorem 7]{Ermakova}, it is enough to show that for any $x,z\in Q^n$, there exists a point $y\in Q^n$ such that the lines $(xy)$ and $(yz)$ belong to $Q^n$. Intersecting with $n-2$ hyperplanes, we can reduce to $n=2$, in which case $Q^2\simeq\P^1\times\P^1$ is covered by two family of lines corresponding to the two rulings. Hence, the point $y=(pr_1(x),pr_2(z))$ satisfies our requirement.}

\demode{Lemma \ref{lem-p2bundle}.}{Since $X$ is smooth and $X\to Y$ is a projective bundle, $Y$ is smooth as well. Since $X$ is rationally connected, $Y$ is rationally connected and by Lemma \ref{fibproj}, one has $-K_Y\cdot C\ge 3$ for any rational curve $C$ in $Y$. By \cite[Cor.1.4]{DeHoe}, $Y$ is a quadric hypersurface $Q^3$ or the projective space $\P^3$. In either case, $Y$ is rational and so it has trivial Brauer group. Hence, $X=\P(V)$ for some vector bundle $V$ on $Y$.

Let us first assume that $Y$ is a quadric hypersurface $Q^3$. Every line $\Delta$ in $Y$ satisfies $-K_Y\cdot \Delta = 3$. Since $-K_X$ has degree at least three on any rational curve, by Lemma \ref{fibproj} and by its equality case, we have $\P(V|_{\Delta})\simeq\P^2\times \Delta$. Hence, for every line $\Delta$ in $Y$, there is an integer $\delta$ such that $V|_{\Delta}$ is isomorphic to $\O_{\P^1}(\delta)^{\oplus 3}$ as a vector bundle on $\Delta\simeq\P^1$. Fixing a single line $\Delta_0$ in $Y$, and noting that $\rho(Y)=1$, we have 
$$3\delta - 3\delta_0 = c_1(V)\cdot \Delta - c_1(V)\cdot \Delta_0 = 0,$$
so the twist $V_0 = V\otimes \O_{Y}(-\delta_0)$ satisfies $V_0|_{\Delta}={\O_{\Delta}}^{\oplus 3}$ for any line $\Delta$ in $Y$. By Lemma \ref{lem-trivquadric}, this vector bundle $V_0$ is trivial on $Y$, and thus $X\simeq \P(V)\simeq\P(V_0)\simeq \P^2\times Q^3$.

Suppose now that $Y$ is a projective space. By the almost-equality case in Lemma \ref{fibproj}, for every line $\Delta$ in $Y$, $$V|_{\Delta}\simeq \bigoplus_{i=1}^3\O_{\P^1}(a_{i,\Delta}),$$
with either $a_{1,\Delta}=a_{2,\Delta}= a_{3,\Delta}$ 
or $a_{1,\Delta}=a_{2,\Delta}= a_{3,\Delta}-1$. Note that the sum $a_{1,\Delta}+a_{2,\Delta}+a_{3,\Delta}=c_1(V)\cdot \Delta$ is independent of the chosen line $\Delta$. If it is divisible by 3, then we are in the first case, else it is congruent to 1 modulo 3 and we are in the second case. In both cases, the $a_{i,\Delta}$ are thus independent of the line $\Delta$.
We fix a line $\Delta_0$ in $\P^3$. The twisted bundle $V_0 = V\otimes {\O_{\P^3}(-a_{1,\Delta_0})}$ now is a uniform bundle of type $(0,0,0)$ or $(0,0,1)$. In the first case, the bundle $V_0$ is globally trivial by \cite{Oko}, and so $X\simeq \P^2\times\P^3$. In the second case, by \cite{EiichiSato}, the vector bundle $V_0$ is either $\O_{\P^3}\oplus\O_{\P^3}\oplus\O_{\P^3}(1)$ or $T_{\P^3}(-1)$, which concludes the classification.}

Let us now study a more general fibred Mori contraction of $X$.

\lem[lem-dimY4]{Let $X$ be a smooth projective rationally connected fivefold and such that, for any rational curve $C\subset X$, one has $-K_X\cdot C \ge 3$. Suppose that $X$ has a fibred Mori contraction $\pi:X\to Y$. Then $\dim Y\le 3$.}

\demo{If $\dim(Y)=4$, the general fibre of $\pi$ is a smooth curve $C$ with trivial normal bundle. By assumption, 
$$2=-K_X\cdot C =\mathrm{deg}_C(-K_C)\ge 3,$$ absurd.}

Let us cover the case when $X$ has a fibred Mori contraction $\pi:X\to Y$ with $1\le \dim(Y)\le 2$.

\lem[lem-dimY2]{Let $X$ be a smooth projective rationally connected fivefold and such that, for any rational curve $C\subset X$, one has $-K_X\cdot C \ge 3$. Suppose that $X$ has a fibred Mori contraction $\pi:X\to Y$ with $1\le \dim Y\le 2$. Then there is a fibred Mori contraction $p:X\to Z$ that is a $\P^2$-bundle.}

\demo{We dispose of a rational curve $C$ such that $\pi(C)\ne\{\mathrm{pt}\}$ and $-K_X\cdot C\ge 3$ is minimal with this condition.
Let ${\cal V}$ be the corresponding family of deformations. By Lemma \ref{lem-unsplit}, ${\cal V}$ is unsplit.
Fix $x\in \mathrm{Locus}({\cal V})$ general. By \cite[Proposition IV.2.6]{Koll} and by assumption, we derive
$$\dim \mathrm{Locus}({\cal V})+ \dim\mathrm{Locus}({\cal V}_x)
\ge -K_X\cdot C + 5 -1
\ge 7.$$
So $\dim \mathrm{Locus}({\cal V}_x) \ge 2.$
By Lemma \ref{lem-finite}, $\dim \mathrm{Locus}({\cal V}_x)\le \dim Y \le 2$.

As equality holds, ${\cal V}$ is a covering family of rational $1$-cycles with $\dim \mathrm{Locus}({\cal V}_x) = 2\ge 5 -3$, so by \cite[Theorem 2, Proposition 1(i)]{BonCaD}, it admits a geometric quotient $p:X\to Z$, that is a fibred Mori contraction, with a general fibre of dimension $2$. If a fibre $F$ of $p$ has dimension 3 or more, then since $\dim Y\le 2$, $\pi|_F$ cannot be finite. So $\pi$ contracts at least a curve $B$ contained in $F$, which is numerically equivalent to a multiple of $C$ as it lies in a ${\cal V}$-equivalence class \cite[Remark 1]{BonCaD}, a contradiction.

So $p$ is an equidimensional fibred Mori contraction with fibres of dimension $2$, of length $-K_X\cdot C\ge 3$. By \cite[Theorem 1.3]{HoeNov}, the morphism $p$ is a $\P^2$-bundle.}

We are left supposing that $X$ has a fibred Mori contraction $\pi:X\to Y$ with $\dim(Y)=3$ that is not a $\P^2$-bundle. Let us first prove a few generalities about its fibres.

\lem{Let $X$ be a smooth projective $n$-dimensional variety with a fibred Mori contraction $\pi$ of length $n-k+1$ onto a variety $Y$ of dimension $k$.
Then the general fibre is isomorphic to $\P^{n-k}$.}

\demo{The general fibre is a smooth variety $F$ of dimension $n-k$ such that $-K_F\cdot C\ge n-k+1$ for any rational curve $C$ in $F$, and $-K_F$ is ample. By \cite{CMSB,Kebekus}, \cite[Theorem 2.1]{HoeNov}, we obtain $F\simeq\P^{n-k}$.}

We recall and prove a fact mentioned in \cite[1.C]{HoeNov}.

\lem[lem-normP3]{Let $X$ be a smooth projective variety of dimension $n\ge 4$ with a fibred Mori contraction $\pi$ of length $n-2$ onto a threefold $Y$. Suppose that $\pi$ is not equidimensional.
Then for any irreducible component $F$ of a fibre of $\pi$ of dimension $n-2$, the normalization $\tilde{F}$ of $F$ is isomorphic to $\P^{n-2}$.}

\demo{By the proof of \cite[Theorem 1.3]{HoeNov}, and as $\mathrm{Univ}_{n-3}(X/Y)\to\mathrm{Chow}_{n-3}(X/Y)$ is a universal family for the $(n-3)$-cycles of $X$ over $Y$, there is a commutative diagram:
$$\xymatrix{
X' \ar[d]_*{\pi'} \ar@/^1pc/[rr]^*{\mu'} \ar[r]_*{\eta'} & \overline{X} \ar[d]_*{\overline{\pi}} \ar[r]_*{\eps'} & X \ar[d]^*{\pi}\\
Y' \ar@/_1pc/[rr]_*{\mu} \ar[r]^*{\eta} & \overline{Y} \ar[r]^*{\eps} & Y
} 
$$

\noindent where $\overline{Y}$ is the normalization of the closure of the $\pi$-equidimensional locus of $Y$ in $\mathrm{Chow}_{n-3}(X/Y)$, $\overline{X}$ is the normalization of the universal family over it, $\eps'$ is the evaluation map, $Y'$ is a resolution of $\overline{Y}$, $X'$ is the corresponding normalized fibred product, $\pi'$ is a $\P^{n-3}$ bundle.
Note that since $Y$ is $\Q$-factorial, the exceptional loci of $\mu$ and of $\eps$ are unions of surfaces, hence the exceptional locus of $\mu'$ is a union of $\P^{n-3}$-bundles on surfaces. Also note that $\pi$, as a fibred Mori contraction, does not contract any divisor; hence, the indeterminacy locus of $\eps^{-1}$ in $Y$ has dimension zero.

Let $F$ be an irreducible component of dimension $n-2$ of a fibre of $\pi$, let $\nu:\tilde{F}\to F$ be its normalization. Let $\Sigma\subset\overline{Y}$ be one of the surfaces that $\eps$ contracts onto $\pi(F)$, chosen such that $\Gamma:=\overline{\pi}^{-1}(\Sigma)$ dominates $F$. Let $S$ be the strict transform of $\Sigma$ by $\eta$, and let $P:=\pi'^{-1}(S)$: it is a $\P^{n-3}$-bundle over $S$ and it dominates $\Gamma$. By the universal property of the normalization, we have a map $f:P\to\tilde{F}$, that fits into the following commutative diagram.

$$\xymatrix{
& & \tilde{F}\ar[d]^*{\nu}\\
 P \ar@/^1pc/[urr]^*{f}\ar[d]_*{\pi'} \ar@/^1pc/[rr]^*{\mu'} \ar[r]_*{\eta'} & \Gamma \ar[d]_*{\overline{\pi}} \ar[r]_*{\eps'} & F \ar[d]^*{\pi}\\
 S \ar@/_1pc/[rr]_*{\mu} \ar[r]^*{\eta} & \Sigma \ar[r]^*{\eps} & \{\mathrm{pt}\}
} 
$$

Let $\ell$ be a line contained in a fibre of $\pi'|_P$. Let ${\cal V}$ be the family of deformation of $f_*\ell$ in $\tilde{F}$.

Let us show that this family satisfies the hypotheses of \cite[Theorem 2.1]{HoeNov}. First, note that $\nu^*(-K_X|_F)$ is ample. Since there is a line in $X'$ numerically equivalent to $\ell$ that is disjoint from all exceptional divisors of $\mu'$, and since $\ell$ is contracted by $\pi'$,
$$\nu^*(-K_X|_F)\cdot f_*\ell=-K_X\cdot \mu'_*\ell=-K_{X'}\cdot\ell=-K_{X'/Y'}\cdot\ell = -K_{\P^{n-3}}\cdot\ell = n-2.$$

Since for any rational curve $C$ in $\tilde{F}$, it holds $\nu^*(-K_X|_F)\cdot C\ge n-2$ by assumption, the family ${\cal V}$ is unsplit. Moreover, it is a covering family, as $\nu$ is birational, $\mu'$ is surjective and the family of deformations of ${\ell}$ is covering. Hence, by \cite[Proposition IV.2.5]{Koll}, for a general point $x\in \tilde{F}$,
$$\dim {\cal V} = n-2 + \dim\mathrm{Locus}{\left({\cal V}_{x}\right)} + 1-3,$$
so we are left to show that $\dim\mathrm{Locus}\left({\cal V}_{x}\right)=n-2$ to conclude.

Let us take $x$ and $y$ general in $F$. It suffices to show that the image by $\mu'|_P$ of a certain fibre $\P^{n-3}$ of $\pi'|_P$ contains both $x$ and $y$, since then there is a line through any two points in $\P^{n-3}$. 

Since $x$ is general and $\Gamma$ dominates $F$, it holds $\dim \eps'^{-1}(x) = \dim \Gamma -\dim F = n-3+2-(n-2)=1$, so there is a one-dimensional family of cycles passing through $x$, parametrized by a curve in $\Sigma$. As there is a finite map $\Sigma\to\mathrm{Chow}_{n-3}(F)$ (a composition of inclusions and a normalization), this is a non-trivial family of divisors. Hence, it must cover $F$, in particular there is one divisor passing through $y$ and $x$. This divisor is dominated by a fibre of $\pi'|_P$, which concludes.}

We now use the fact that $\pi$ is not a $\P^2$-bundle (in fact, that $\pi$ is not equidimensional) to construct covering families of rational curves on $X$. Before that, we prove a simple lemma.

\defi{Let $f:X\dashrightarrow Y$ be a rational map. We say that $f$ is {\it almost holomorphic} if there are Zariski open subsets $U\subset X$ and $V\subset Y$ such that $f|_U: U\to V$ is a proper holomorphic map.}

\lem[lem-almostholo]{Let $f:X\dashrightarrow Y$ be almost holomorphic map. If $Y$ is a curve, then $f$ is holomorphic.}

\demo{Let $\eps:X'\to X$ be a resolution of indeterminacies for $f$, let $f':X'\to Y$ be the induced holomorphic map. As $f$ is almost holomorphic, no component of the exceptional locus of $\eps$ is dominant onto $Y$. As $Y$ is curve, this means that the exceptional locus of $\eps$ is sent onto finitely many points in $Y$. So $f'$ factors through $\eps$, i.e., $f$ is holomorphic.}

\lem[lem-covers]{Let $X$ be a smooth projective rationally connected fivefold, such that $-K_X\cdot C\ge 3$ for any rational curve $C\subset X$. Suppose that $X$ has a fibred Mori contraction $\pi:X\to Y$ with $\dim Y=3$. If $\pi$ is not a $\P^2$-bundle, then any rational curve $C\subset X$ such that $\pi(C) \ne \{\mathrm{pt}\}$, and which deforms in an unsplit family, deforms in a family covering $X$.}

\demo{Note that if $\pi$ is equidimensional, by \cite[Theorem 1.3]{HoeNov} it is a $\P^2$-bundle. Hence, we assume that a variety $F$ of dimension 3 is contained in a fibre of $\pi$. By contradiction, we consider a rational curve $C\subset X$ such that $\pi(C) \ne \{\mathrm{pt}\}$, and the family ${\cal V}$ of deformations of $C$ is unsplit and not covering $X$.

Fix $x\in\mathrm{Locus}({\cal V})$ general. By Lemma \ref{lem-finite}, $\dim \mathrm{Locus}({\cal V}_x)\le \dim Y \le 3$.
Since the family ${\cal V}$ is unsplit, 
$$\dim \mathrm{Locus}({\cal V})+ \dim\mathrm{Locus}({\cal V}_x)
\ge -K_X\cdot C + 5 -1
\ge 7,$$
in particular as ${\cal V}$ is not covering, $\dim\mathrm{Locus}({\cal V})=4$ and $\dim\mathrm{Locus}({\cal V}_x)=3$. 

Let $n:\tilde{D}\to D$ denote the normalization of $D=\mathrm{Locus}({\cal V})$, and let $\tilde{{\cal V}}$ be the covering family on $\tilde{D}$. Note that $\pi$ induces a fibration of $\tilde{D}$ onto a variety of smaller dimension that is not a point, in particular $\rho(\tilde{D})\ge 2$. Thus, by \cite[Corollary 4.4]{ACOMuk}, $\tilde{D}$ cannot be $\tilde{{\cal V}}$-chain-connected. 

Considering the dominant almost holomorphic map $r:\tilde{D}\dashrightarrow Z$ whose general fibre is a $\tilde{{\cal V}}$-equivalence class \cite[Section 2]{BonCaD}, the variety $Z$ is thus not a point. Since $\dim\mathrm{Locus}(\tilde{{\cal V}}_x)=3$ for a general $x\in\mathrm{Locus}(\tilde{{\cal V}})$, the variety $Z$ must be a curve, in particular, by Lemma \ref{lem-almostholo}, the map $r$ is holomorphic.

Note that, as $D$ is a relatively ample Cartier divisor with respect to $\pi$, it intersects the three-dimensional variety $F$ along a surface $S$. Since $\dim n^{-1}(S)=2 >\dim Z =1$, the restriction $r|_{n^{-1}(S)}:n^{-1}(S)\to Z$ cannot be finite. So it contracts a curve $B$. Its image $n(B)$ is in a ${\cal V}$-equivalence class, so as ${\cal V}$ is unsplit, it is numerically equivalent to a multiple of $C$. But $n(B)\subset F$, so this curve is contracted by $\pi$, a contradiction.}

\defi{Let $f:X\to Y$ be a morphism of normal varieties. We say that $f$ is {\it quasi{\'e}tale} if $\dim X=\dim Y$, and there is a Zariski closed subset $Z$ in $X$ of codimension at least 2 such that $f:X\setminus Z\to Y\setminus f(Z)$ is {\'e}tale.}

\rem{Note that if $f:X\to Y$ is a finite quasi{\'e}tale cover and $Y$ is smooth, then by Zariski purity of the branch locus \cite[Proposition 2]{ZarPur}, $f$ is {\'e}tale.}

\lem[lem-P3]{Let $X$ be a smooth projective rationally connected fivefold, such that $-K_X\cdot C\ge 3$ for any rational curve $C\subset X$. Suppose that $X$ has a fibred Mori contraction $\pi:X\to Y$ with $\dim Y > 0$. If $X$ is not a $\P^2$-bundle over any smooth projective base, then $Y\simeq\P^3$. Moreover, we have $\rho(X)=2$, and if $C$ is a line in the smooth locus of $\pi$ in $Y$ and $s$ is a minimal section over $C$ in $X$, the class of $s(\P^1)$ generates the other extremal ray in $\overline{NE}(X)$, induces a fibred Mori contraction to a positive dimensional variety too, and satisfies $-K_X\cdot s(\P^1) = 3$.}

\demo{Note that $\dim(Y) = 3$, by Lemmas \ref{lem-dimY4}, \ref{lem-dimY2}. By the last lemma of \cite{DruelErratum}, let $C$ be a minimal free rational curve in the smooth locus $Y^0\subset Y$ of $\pi$. Let $s$ be a minimal section over $C$. Lemma \ref{fibproj} yields 
$$4\ge -K_Y\cdot C\ge -K_X\cdot s(\P^1).$$
The family ${\cal V}$ of deformations of $s(\P^1)$ is unsplit.
Indeed, suppose by contradiction that it is splitting. By \cite[Explanation IV.2.2]{Koll} there is a cycle
$$\sum_{i} a_iC_i\eqnum s(\P^1),$$
with $C_i$ rational curves, $a_i\ge 1$ integers, and $\sum_i a_i\ge 2$.
Then, intersecting with $-K_X$ yields $4\ge -K_X\cdot s(\P^1)\ge 6$, a contradiction.

 By Lemma \ref{lem-covers}, ${\cal V}$ therefore is a covering family. By \cite[Proposition IV.2.6]{Koll}, it moreover holds 
$$\dim \mathrm{Locus}({\cal V}_x) \ge -K_X\cdot s(\P^1) - 1\ge 2 = 5 -3,$$
so by \cite[Theorem 2, Proposition 1(i)]{BonCaD}, there is a geometric quotient $p:X\to Z$, that is a fibred Mori contraction, with general fibre of dimension at least $-K_X\cdot s(\P^1) - 1$. By Lemma \ref{lem-dimY4}, we have $\dim Z\le 3$ and by Lemma \ref{lem-dimY2}, we have $\dim(Z) = 3$, or $X$ is a $\P^2$-bundle over some three-dimensional base. So $\dim Z = 3$, hence $-K_X\cdot s(\P^1)= 3$. It also follows that $s(\P^1)$ is an extremal class in the Mori cone, as wished.

Again, $X$ not being a $\P^2$-bundle over any smooth base, $p$ is not equidimensional by \cite[Theorem 1.3]{HoeNov}, so a variety $F$ of dimension 3 is contained in a fibre of $p$. By Lemma \ref{lem-normP3}, the normalization $n:\tilde{F}\to F$ satisfies $\tilde{F}\simeq\P^3$.

Since $\pi$ and $p$ are distinct Mori contractions, they contract no common numerical class of curve, in particular $\pi|_F:F\to Y$ is finite onto its image, hence finite surjective for dimensional reasons. There is an effective ramification divisor $R\in\mathrm{Pic}(\P^3)$ such that
$-K_{\P^3}= n^*{\pi|_F}^*(-K_Y)-R$. As $F$ is an irreducible component of a ${\cal V}$-equivalence class, and as ${\cal V}$ is unsplit, $F$ contains a deformation of $s(\P^1)$. Let $\tilde{C}$ be the lift to $\tilde{F}$ of a deformation of $s(\P^1)$ that is contained in $F$. Then $-K_{\P^3}\cdot\tilde{C}\ge 4$, and $n^*{\pi|_F}^*(-K_Y)\cdot\tilde{C}=-K_Y\cdot C\le 4$. So $R\cdot\tilde{C}\le 0$, but $R\in\mathrm{Pic}(\P^3)$ is effective, thus ample or trivial, so $R$ is trivial. The finite map $\pi|_{F}\circ n :\P^3\to Y$ is thus quasi{\'e}tale. So, its base change $\P^3\underset{Y}{\times} X\to X$ is also quasi{\'e}tale, as $\pi:X\to Y$ contracts no divisor. But $X$ is rationally connected, hence simply connected, and smooth, so $\P^3\underset{Y}{\times} X\to X$ is an isomorphism. Hence $\pi|_{F}\circ n:\P^3\to Y$ is an isomorphism too. 

Since $\rho(Y) = 1$, we have $\rho(X) = 2$. 
Since $Y\simeq \P^3$ and $4\ge -K_Y\cdot C$, the curve $C$ is a line.}

\lem[lem-twoP3]{Let $X$ be a smooth projective rationally connected fivefold, such that $-K_X\cdot C\ge 3$ for any rational curve $C\subset X$. Suppose that $X$ has a fibred Mori contraction $\pi:X\to Y$ with $\dim(Y)>0$. If $X$ is not a $\P^2$-bundle over any smooth projective base, then $\rho(X)=2$ and $X$ has two distinct fibred Mori contractions onto $\P^3$, with corresponding extremal rays generated by the minimal sections $s(\P^1),\sigma(\P^1)$ above lines that lie in each $\P^3$ in the smooth locus of the fibration. Moreover, 
$$-K_X\cdot s(\P^1)=-K_X\cdot \sigma(\P^1)=3.$$}

\demo{Apply Lemma \ref{lem-P3} twice.}

\demode{Proposition \ref{prop-summary}}{If $X$ has a $\P^2$-bundle structure, then Lemma \ref{lem-p2bundle} concludes. Suppose that $X$ is not a $\P^2$-bundle. By Lemma \ref{lem-twoP3}, $X$ admits exactly two fibred Mori contractions $\pi$ and $p$, both onto $\P^3$. 
Given the intersection number of $-K_X$ with both extremal rays, and as $\pi_*s(\P^1)$ is a line in $\P^3$ and as $p_* s(\P^1)=0$, we have
$$-K_X\cdot s(\P^1)=3=\pi^*\O_{\P^3}(3)\cdot s(\P^1)=(\pi^*\O_{\P^3}(3)\otimes p^*\O_{\P^3}(3))\cdot s(\P^1),$$
and similarly $$-K_X\cdot \sigma(\P^1)=(\pi^*\O_{\P^3}(3)\otimes p^*\O_{\P^3}(3))\cdot \sigma(\P^1).$$
Hence, as $\rho(X)=2$, and $s(\P^1)$ and $\sigma(\P^1)$ are independent,
$${\omega_X}^* = \pi^*\O_{\P^3}(3)\otimes p^*\O_{\P^3}(3).$$ By Theorem \ref{theo-n-1main}, $-K_X$ is ample. So $X$ is a Fano fivefold, and we just showed that it has index 3. By the classification in \cite{WisIndex}, $X$ must then be a $\P^2$-bundle, which is a contradiction.}

\subsection{Divisorial contractions}

Let us classify divisorial Mori contraction of large length.

\prop[prop-nobiglength]{Let $X$ be a smooth projective rationally connected variety of dimension $n$ such that $-K_X\cdot C\ge 3$ for every rational curve $C$. Then $X$ admits no divisorial Mori contraction of length greater or equal to $n-1$.}

\rem{In particular, the assumptions are fulfilled if there is $1\le r\le n-1$ such that $\bigwedge^r T_X$ is strictly nef, by \cite[Theorem 1.2]{Main} and Lemma \ref{degrantican}.}

The proof uses the following lemma, that excludes some special contractions of length $n-1$.

\lem[lem-noblowup]{Let $X$ be a smooth projective rationally connected variety of dimension $n$ such that $-K_X\cdot C\ge 3$ for every rational curve $C$. Then there is no morphism $X\to Y$ that is a blow-up of a smooth point in a smooth variety.}

\demode{Lemma \ref{lem-noblowup}}{By contradiction, consider such a smooth blow-up:
$$f: E\subset X\to p\in Y$$

Note that since $X$ is rationally connected, so $Y$ is too.
Let $C$ be a rational curve through $p$. 

Since $-f^*K_Y=-K_X+(n-1)E$ and since no curve is contained in the blown-up locus $p$, the anticanonical divisor $-K_Y$ is stricly nef. By bend-and-break \cite[Proposition 3.2]{Debarre} on the smooth variety $Y$, one can thus assume $-K_Y\cdot C\le n+1$. The strict transform $C'\subset X$ of $C$ satisfies $E\cdot C' > 0$.
Since $K_X=f^*K_Y+(n -1)E$, we have
 $$3\le -K_X\cdot C'\le -K_Y\cdot C - (n-1)\le 2,$$
a contradiction!}

\demode{Proposition \ref{prop-nobiglength}}{By Ionescu-Wi\'snewski inequality, if $X$ admits a divisorial Mori contraction of length $\ell\ge n-1$, the exceptional divisor $E$ and the general fibre $F\subset E$ satisfy:
$$\dim E+\dim F \ge n +\ell - 1\ge 2n-2,$$ i.e., $\ell=n-1$ and $E=F$ is contracted onto a point.
So \cite[Theorem 5.2]{AO} applies and shows that this divisorial Mori contraction of $X$ correponds to a blow-up of a smooth point in a smooth variety, which contradicts Lemma \ref{lem-noblowup}.}

We now consider divisorial Mori contractions of length $n-2$.

\prop[prop-nodivtopt]{Let $X$ be a smooth projective variety of dimension $n\ge 5$, that is rationally connected and such that $-K_X\cdot C\ge n-2$ for any rational curve $C\subset X$. Then $X$ has no divisorial Mori contraction contracting the exceptional divisor to a point.}

\rem{These assumptions are fulfilled if $\bigwedge^4 T_X$ is strictly nef, by \cite[Theorem 1.2]{Main} and Lemma \ref{degrantican}.}

\demo{Assume that $\eps: X\to Y$ is a divisorial Mori contraction contracting the exceptional divisor $E$ to a point. Note that as $X$ is rationally connected, there exists a rational curve $C$ that intersects $E$ without being contained in $E$. In particular, $E\cdot C > 0$. Among all such curves, let actually $C$ be one such that $-K_X\cdot C$ is minimal. Then we claim that the family ${\cal V}$ of deformations of $C$ is unsplit.
Indeed, suppose by contradiction that it is splitting. By \cite[Explanation IV.2.2]{Koll}, we have
$$C\eqnum\sum_i a_iC_i,$$
with rational curves $C_i$ and coefficients $a_i\ge 1$ such that $\sum a_i\ge 2$.
Then $E\cdot C > 0$, so without loss of generality, $E\cdot C_1 > 0$. In particular, $C_1$ intersects $E$ and is not contracted by $\eps$, hence not contained in $E$. Since $-K_X$ has positive degree on all rational curves in $X$, we have $-K_X\cdot C_1 < -K_X\cdot C$, which contradicts the minimality of $-K_X\cdot C$.

By \cite[Proposition IV.2.6.1]{Koll}, for a general $x\in\mathrm{Locus}({\cal V})$,
 $$\dim\mathrm{Locus}({\cal V})+\dim\mathrm{Locus}({\cal V}_x)\ge n+n-2-1.$$
 
In particular, $\dim\mathrm{Locus}({\cal V}_x)\ge n-3$, and as $X$ is smooth, $E$ is Cartier, hence intersects $\mathrm{Locus}({\cal V}_x)$ along a subscheme of dimension at least $n-4\ge 1$. Let $B$ be a curve in this intersection. It is contained in $E$, hence contracted by $\eps$, hence satisfies $E\cdot B<0$. On the other hand, it is contained in $\mathrm{Locus}({\cal V}_x)$, hence is numerically equivalent to a multiple of $C$ by \cite[Lemma 4.1]{ACOMuk}. It has to be a positive multiple, as one sees when intersecting with any ample divisor. But $E\cdot C > 0$, a contradiction.}

\cor[cor-blp5]{Let $X$ be a smooth projective variety of dimension $n\ge 5$, that is rationally connected and such that $-K_X\cdot C\ge n-2$ for any rational curve $C\subset X$. Suppose that $\eps:X\to Y$ is a divisorial Mori contraction. Then $Y$ is smooth and $\eps$ is the blow-up of a smooth curve in $Y$.}

\demo{Recall \cite[Proposition 6.10(b)]{Debarre} that the divisorial Mori contraction $\eps$ has a unique exceptional divisor $E$ as its exceptional locus.  By \cite[Lemma 2.62]{KollarMori}, a ray $\R_+[C]$ associated to $\eps$ satisfies $E\cdot C < 0$, so such $C$ has negative intersection with at least one effective divisor. Moreover, $\eps$ is a Mori contraction of length $n-2$. So \cite[Theorem 5.3]{AO} applies, showing that $\eps$ either contracts a divisor to a point, or is a blow-up of a smooth curve in a smooth variety $Y$. By Proposition \ref{prop-nodivtopt}, only the latter can occur.}

Let us finally describe more precisely what happens in the occurrence of Corollary \ref{cor-blp5}.

\lem[lem-r=n-1]{Let $X$ be a smooth projective variety of dimension $n\ge 3$, that is rationally connected and such that for some $1\le r\le n-1$, for any rational curve $C\subset X$, it holds $-K_X\cdot C\ge n+2-r$. If there is a morphism $\eps:X\to Y$ that is a blow-up of a smooth curve in the smooth variety $Y$, then $r=n-1$.}

\demo{Consider such a smooth blow-up:
$$f: E\subset X\to \ell\subset Y$$

As $X$ is rationally connected, so is $Y$. Fix $H$ an ample divisor on $Y$.
Let $C\subset Y$ be a rational curve other than $\ell$ passing through a point $p\in\ell$, with $H\cdot C$ minimal among the degrees of all rational curves intersecting $\ell$ other than $\ell$. Fix another point $q\in C\setminus C\cap\ell$. By bend-and-break \cite[Proposition 7.3]{Debarre}, as $Y$ is smooth, if $-K_Y\cdot C\ge n+2$, then there is a connected non-integral 1-cycle that is a deformation of $C$ passing through $p$ and $q$. In particular, $$\sum_{i=1}^k a_iC_i\eqnum C,$$
with rational curves $C_i$ such that $p\in C_1,\, q\in C_{i_0}$ for some $i_0$, coefficients $a_i\ge 1$, and $\sum_{i=1}^k a_i\ge 2$. As $q\not\in\ell$, we have that $C_{i_0}\ne\ell$, so either $C_1\ne\ell$, or $C_1=\ell$ and $k\ge 2$.
Intersecting with $H$, we see that $H\cdot C_i < H\cdot C$ for all $i$, in particular for $C_1$. If $C_1\ne \ell$, then $H\cdot C_1$ contradicts the minimality of $H\cdot C$. If $C_1=\ell$, then $k\ge 2$ and by connectedness of the rational cycle, there is a curve $C_{i_1}\ne \ell$ that intersects $C_1=\ell$. So $C_{i_1}\ne\ell$ intersects $\ell$ and contradicts the minimality, as $H\cdot C_{i_1}<H\cdot C$ again.
So $-K_Y\cdot C\le n+1$.

 The strict transform $C'\subset X$ of $C$ satisfies $E\cdot C' > 0$.
Since $K_X=f^*K_Y+(n -2)E$, and by assumption,
 $$n+2-r\le -K_X\cdot C'\le -K_Y\cdot C - (n-2)\le 3,$$
 so $r=n-1$.}

\prop[prop-fibred]{Let $X$ be a smooth projective variety of dimension $n\ge 5$, that is rationally connected and such that $\bigwedge^4 T_X $ is strictly nef. If there is a morphism $\eps:X\to Y$ that is a blow-up of a smooth curve in the smooth variety $Y$, then $X$ is a fivefold and
there is a fibred Mori contraction $\pi: X\to Z$ with $\dim(Z)>0$.}

\demo{By Lemma \ref{lem-r=n-1}, we have $n=5$. So by Theorem \ref{theo-n-1main}, $-K_X$ is ample. The Mori cone $NE(X)$ is closed, generated by finitely many classes of rational curves. Let $E$ be the exceptional divisor of $\eps$. Note that there exists an extremal ray $R=\R_+[C]$ of $NE(X)$ on which $E\cdot C > 0$. Indeed, if there were not such a ray, then $E$ would be non-positive on all curves in $X$, which is absurd for an effective divisor. So, let $R=\R_+[C]$ be an extremal ray on which $E\cdot C > 0$.

Denote the associated Mori contraction by $\pi : X\to Z$. Since $X$ already had a non-trivial Mori contraction $\eps$, we have $\dim(Z)>0$. Let us prove that $\pi$ is a fibred Mori contraction.

By Lemma \ref{paspet}, $\pi$ cannot be a small contraction. Assume by contradiction that it is a divisorial contraction. By Corollary \ref{cor-blp5}, the variety $Z$ is smooth and $\pi$ is a blow-up along a smooth curve of $Z$. Let $E'$ be the $\pi$-exceptional divisor. Let $\ell$, respectively $\ell'$, be the image of $E$, respectively $E'$, in $Y$, respectively $Z$. Let $F'$ be a general fibre of $\pi|_{E'}$. It has dimension $n-2$.
Note that $F'$ and $E$ intersect, since $E\cdot C > 0$.
Hence, $E\cap F'$ is a subscheme of $X$ of dimension at least $n-3$. Since $\eps$ and $\pi$ are distinct Mori contractions, the restriction $\eps|_{E\cap F'}$ must be finite onto its image, which is contained in $\ell$. So $n-3\le 1$, a contradiction!

So $\pi$ is a fibred Mori contraction.}

\prop[prop-blpn]{Let $X$ be a smooth projective variety of dimension $n\ge 5$, that is rationally connected and such that $\bigwedge^4 T_X $ is strictly nef. If there is a morphism $\eps:X\to Y$ that is a blow-up of a smooth curve, then $Y\simeq\P^5$ and $\eps$ is the blow-up of a line.}

\demo{By Proposition \ref{prop-fibred}, $X$ is a fivefold and admits a fibred Mori contraction onto a positive dimensional base. So Proposition \ref{prop-summary} applies, showing that $X$ belongs to a list of certain varieties of Picard number two. Only one of them has a divisorial Mori contraction, namely $\mathrm{Bl}_{\ell}(\P^5)=\P(\O_{\P^3}\oplus \O_{\P^3}\oplus \O_{\P^3}(1))$.}

\section{Results on $\bigwedge ^3 T_X$}\label{sec-3}

\demode{Theorem \ref{theo-3mainnew}}{Note that $-K_X$ is nef, and non-trivial (as it is positive on rational curves, and $X$ is rationally connected). If $\rho(X)=1$, $-K_X$ is ample and $X$ is thus a Fano variety. If $\rho(X)\ge 2$, by the Cone Theorem, $X$ admits a Mori contraction, which by Lemma \ref{paspet} and Proposition \ref{prop-nobiglength} is a fibred Mori contraction. Corollary \ref{cor-evid} implies that $X$ is a fourfold. By Lemma \ref{prop-p2p2}, $X$ has an equidimensional fibred Mori contraction to a surface, so by Proposition \ref{prop-genfib}, we have $X\simeq\P^2\times\P^2$.}

\demode{Corollary \ref{cor-3main}}{It is straightforward from Lemma \ref{degrantican}, \cite[Theorem 1.2]{Main}, and Theorem \ref{theo-3mainnew}.}

\rem{It is easy to check that $\bigwedge^3 T_{\P^2\times\P^2}$ is in fact ample.}

\expl{Let $X$ be a cubic in $\P^n$ with $n\ge 5$. From the tangent exact sequence
$$0\to T_X\to T_{\P^n}|_X\to \O_X(3)\to 0,$$
we can use \cite[II.Ex.5.16(d)]{HarBook} to derive the existence of a surjection
$$0\to F_4\to \bigwedge^4 T_{\P^n}|_X\to \bigwedge^3 T_X\otimes\O_X(3)\to 0.$$
 As $T_{\P^n}|_X\otimes\O_X(-1)$ is nef, the quotient of its fourth exterior power $\bigwedge^3 T_X\otimes\O_X(-1)$ is also nef, and thus $\bigwedge^3 T_X$ is ample.}

\expl{Let $X$ be the complete intersection of two quadrics in $\P^n$ with $n\ge 6$. From the tangent exact sequence
$$0\to T_X\to T_{\P^n}|_X\to \O_X(2)\oplus\O_X(2)\to 0,$$
we can use \cite[II.Ex.5.16(d)]{HarBook} to derive the existence of a surjection
$$0\to F_4\to \bigwedge^5 T_{\P^n}|_X\to \bigwedge^3 T_X\otimes\O_X(4)\to 0.$$
As $T_{\P^n}|_X\otimes\O_X(-1)$ is nef, the quotient of its fifth exterior power $\bigwedge^3 T_X\otimes\O_X(-1)$ is also nef, and thus $\bigwedge^3 T_X$ is ample.}

\section{Results on $\bigwedge^4 T_X$}\label{sec-4}

\subsection{Examples}

\lem{Let $X$ be the fivefold $\P(T_{\P^3})$. Then $\bigwedge^4T_X$ is ample.}

\demo{Denote the natural projection by $p:X\to \P^3$, the tautological line bundle on $X$ by $\O_X(1)$. By \cite[II.Ex.5.16(d)]{HarBook}, there is an exact sequence
$$0\to \bigwedge^2 T_{X/\P^3}\otimes p^*\bigwedge^2 T_{\P^3}
\to \bigwedge^4 T_X
\to T_{X/\P^3}\otimes p^*\O_{\P^3}(-K_{\P^3})
\to 0.$$

 Let us prove that $E_1 = T_{X/\P^3}\otimes p^*\O_{\P^3}(-K_{\P^3})$ is ample. 
We have the relative Euler sequence
$$0\to \O_X\to p^*\Omega^1_{\P^3}\otimes \O_X(1)\to T_{X/\P^3}\to 0.$$
The bundle $E_1$ is a quotient of $p^*\Omega^1_{\P^3}(4)\otimes \O_X(1)$. But as $T_{\P^3}$ is ample, $\O_X(1)$ is ample. Moreover, $\Omega^1_{\P^3}(4)\simeq \bigwedge^2 T_{\P^3}$ is ample too, which concludes by \cite[6.1.16]{Laz}.

Let us prove that $E_2=\bigwedge^2 T_{X/\P^3}\otimes p^*\bigwedge^2 T_{\P^3}$ is ample. This would settle the ampleness of $\bigwedge^4 T_X$ by \cite[6.1.13(ii)]{Laz}.
From \cite[II.Ex.5.16(d)]{HarBook} and the relative Euler sequence, we derive 
$$0\to T_{X/\P^3}\to p^*T_{\P^3}(-4)\otimes \O_X(2)\to \bigwedge^2 T_{X/\P^3}\to 0.$$
Since $E_2$ is a quotient of $p^*(T_{\P^3}(-4)\otimes\bigwedge^2 T_{\P^3})\otimes \O_X(2)$, we are left proving that the latter is ample. Notice that $T_{\P^3}(-1)$ is globally generated and thus nef. So the bundle $T_{\P^3}(-3)\otimes\bigwedge^2 T_{\P^3}=T_{\P^3}(-1)\otimes\bigwedge^2 T_{\P^3}(-1)$ is nef as well. Finally, $\O_X(1)$ is ample, and we see that $\O_X(1)\otimes p^*\O_{\P^3}(-1)$ is a quotient of $p^*T_{\P^3}(-1)$ (dualizing the relative Euler exact sequence and twisting by $\O_X(1)$), hence it is nef. We conclude by \cite[6.2.12(iv)]{Laz}.}

\lem{Let $X$ be the fivefold $\P(\O_{\P^3}\oplus \O_{\P^3}\oplus \O_{\P^3}(1))$. Then $\bigwedge^4T_X$ is ample.}

\rem{Note that $\P(\O_{\P^3}\oplus \O_{\P^3}\oplus \O_{\P^3}(1))$ is isomorphic to the blow-up of  line in $\P^5$ \cite[Section 9.3.2]{3264}.}

\demo{Denote the natural projection by $p:X\to \P^3$, the tautological line bundle on $X$ by $\O_X(1)$. By \cite[II.Ex.5.16(d)]{HarBook}, there is an exact sequence
$$0\to \bigwedge^2 T_{X/\P^3}\otimes p^*\bigwedge^2 T_{\P^3}
\to \bigwedge^4 T_X
\to T_{X/\P^3}\otimes p^*\O_{\P^3}(-K_{\P^3})
\to 0.$$

Let us prove that $E_1 = T_{X/\P^3}\otimes p^*\O_{\P^3}(-K_{\P^3})$ is ample. 
We have the relative Euler sequence
$$0\to \O_X\to p^*(\O_{\P^3}\oplus \O_{\P^3}\oplus \O_{\P^3}(-1))\otimes \O_X(1)\to T_{X/\P^3}\to 0.$$
The bundle $E_1$ is a quotient of $p^*(\O_{\P^3}(3)\oplus \O_{\P^3}(4)\oplus \O_{\P^3}(4))\otimes \O_X(1)$. Since $\O_{\P^3}(3)\oplus \O_{\P^3}(4)\oplus \O_{\P^3}(4)$ is ample and $\O_X(1)$ is nef and $p$-ample, the bundle $E_1$ is thus ample.

Let us prove that $E_2 = \bigwedge^2 T_{X/\P^3}\otimes p^*\bigwedge^2 T_{\P^3}$ is ample. 
From \cite[II.Ex.5.16(d)]{HarBook} and the relative Euler sequence, we derive 
$$0\to T_{X/\P^3}\to p^*(\O_{\P^3}(-1)\oplus \O_{\P^3}(-1)\oplus\O_{\P^3})\otimes \O_X(2)\to \bigwedge^2 T_{X/\P^3}\to 0.$$

It is thus enough to prove that $p^*\bigwedge^2 T_{\P^3}\otimes p^*\O_{\P^3}(-1)\otimes \O_X(2)$ is ample, which is clear since $\bigwedge^2 T_{\P^3}(-1)=(\bigwedge^2 T_{\P^3})(-2)$ is globally generated and thus nef, and since $p^*\O_{\P^3}(1)\otimes \O_X(2)$ is ample.}

\rem{It is easy check to that $\bigwedge ^4 T_{\P^2\times\P^3}$, $\bigwedge ^4 T_{\P^2\times Q^3}$, $\bigwedge ^4 T_{\P^3\times\P^3}$ are ample.}

\subsection{Proof of Theorem \ref{theo-4mainnew} and Corollary \ref{cor-4main}}

\demode{Theorem \ref{theo-4mainnew}}{Note that $-K_X$ is nef, and non-trivial (as it is positive on rational curves, and $X$ is rationally connected). If $\rho(X)=1$, $-K_X$ is ample and $X$ is thus a Fano variety. If $\rho(X)\ge 2$, by the Cone Theorem, $X$ admits a Mori contraction. By Lemma \ref{paspet}, it cannot be a small contraction.

Suppose that it is a divisorial contraction. By Corollary \ref{cor-blp5} and Lemma \ref{lem-r=n-1}, it is a smooth blow-up of a smooth curve in a fivefold, but we are assuming that $X$ has dimension at least six, a contradiction!

So $X$ has no divisorial contraction. Thus, it has a fibred Mori contraction onto a positive dimensional variety. Corollary \ref{cor-evid} implies that $X$ is a fivefold or a sixfold. 
By assumption, $X$ is thus a sixfold. By Lemma \ref{prop-p2p2}, $X$ has an equidimensional fibred Mori contraction to a threefold, so by Proposition \ref{prop-genfib}, we have $X\simeq\P^3\times\P^3$, which concludes.}

\demode{Corollary \ref{cor-4main}}{By Theorem \ref{theo-4mainnew}, is is enough to consider the case when $X$ is a fivefold. In particular, by Theorem \ref{theo-n-1main}, $X$ is a Fano variety.
Again, if $\rho(X)=1$, there is nothing to prove.

If $\rho(X)\ge 2$, by the Cone Theorem, $X$ admits a Mori contraction. By Lemma \ref{paspet}, it cannot be a small contraction.

Suppose that it is a divisorial contraction. By Corollary \ref{cor-blp5}, it is a smooth blow-up of a smooth curve, and by Proposition \ref{prop-blpn}, $X\simeq\mathrm{Bl}_{\ell}\P^5$.

Otherwise, it is a fibred Mori contraction onto a positive dimensional variety. Since $X$ is a fivefold such that $\bigwedge^4 T_X$ is strictly nef, Proposition \ref{prop-summary} applies and concludes.}

\newcommand{\etalchar}[1]{$^{#1}$}

\end{document}